\documentclass {article}

\usepackage{amsfonts,amsmath,latexsym}
\usepackage{amstext,amssymb}
\usepackage{amsthm}
\usepackage[scanall]{psfrag}
\usepackage {graphicx}
\usepackage{verbatim}

\newcommand {\R}{\mathbb{R}}

\newcommand {\C}{\mathbb{C}}
\newcommand {\D}{\mathbb{D}}

\newcommand {\del}{\partial}

\newtheorem {thm} {Theorem}
\newtheorem*{thmA}{Theorem A}

\newtheorem {lemma} [thm] {Lemma}
\newtheorem {cor} [thm] {Corollary}
\newtheorem {defn} {Definition}
\newtheorem {rmk} {Remark}

\newcommand{\beq}{\begin{equation}}
\newcommand{\eeq}{\end{equation}}

\begin {document}

\title{ 
{\bf Isoperimetric Inequalities and 
Variations on Schwarz's Lemma}}
\author{Tom Carroll 
and Jesse Ratzkin \\
University College Cork and University 
of Cape Town\\
{\tt t.carroll@ucc.ie}
and {\tt j.ratzkin@ucc.ie}} 
\maketitle

\begin {abstract} 
In this note we prove a version of the 
classical Schwarz lemma for the first eigenvalue 
of the Laplacian with Dirichlet boundary data. 
A key ingredient in our proof is an isoperimetric 
inequality for the first eigenfunction, due to Payne 
and Rayner, which we reinterpret as an isoperimetric 
inequality for a (singular) conformal metric 
on a bounded domain in the plane. 
\end {abstract}

\section{Introduction}

Let $\D = \{ z : |z| < 1\} \subset \C$ be the 
unit disk in the complex plane, and let 
$f:\D \rightarrow \D$ be analytic with $f(0) = 0$. 
Then the  
classical Schwarz lemma states that $|f(z)| 
\leq |z|$ and $|f'(0)| \leq 1$, and that equality
in either case implies $f(z) = e^{i\theta} z$ 
for some $\theta \in \R$. One can reinterpret 
this result more geometrically by defining, 
for $0<r<1$, 
$$\operatorname {Rad}(r) = \sup_{|z| = r} 
|f(z) - f(0)|,$$
so that the Schwarz lemma states $\operatorname 
{Rad}(r) \leq r\operatorname{Rad}(1)$ for 
every analytic $f:\D \rightarrow \C$. In fact, the classical 
proof of the Schwarz lemma implies 
$$\Phi_{\operatorname {Rad}} (r) = 
\frac{\operatorname {Rad}(r)}{r} $$
is a strictly increasing function of $r$, 
unless $f$ is linear (in which case 
$\Phi_{\operatorname{Rad}}$ is constant). 

Burckel, Marshall, Minda, Poggi-Corradini, and 
Ransford \cite{BMMPR} recently 
proved versions of the Schwarz Lemma for diameter, 
logarithmic capacity, and area. They asked whether 
similar inequalities hold for other quantities, 
such as the first eigenvalue $\lambda$ of the Laplacian 
with Dirichlet boundary data. 
\begin {thm} \label{eigen-schwarz}
Let $f$ be a conformal mapping of the unit disk 
$\D$. The function
\beq\label{1.1}
\Phi_\lambda(r) =  
\frac{\lambda\big( f(r\D) \big)}{\lambda(r\D)}
 	= \frac{1}{j_0^2} r^2 \lambda\big( f(r\D) \big), 
\quad 0 < r < 1,
\eeq
is strictly decreasing, unless $f$ is linear (in which 
case $\Phi_\lambda$ is constant). 
\end {thm}
Taking a limit of the right hand side of \eqref{1.1} as $r\rightarrow 0^+$, we 
recover the estimate in Section 5.8 of \cite{PS}. 
Using monotonicity of the right 
hand side of \eqref{1.1}, we also see that the limit as
$r \rightarrow 1^-$ of the right hand side of \eqref{1.1} 
also exists, thought it might be zero.

A slight modification of the proof of Theorem~\ref{eigen-schwarz} yields
the following corollary. 
\begin {cor} \label {eigen-schwarz2}
Let $f$ be an analytic function in the unit disk. 
For $0 < r < 1$ let 
$\Sigma_r$ be the Riemann surface associated to $f:r\D \rightarrow 
\C$. Then the function 
$$ r \mapsto  \frac{1}{j_0^2} r^2 \lambda(\Sigma_r)$$
is strictly decreasing, unless $f$ is linear (in 
which case this function is constant). 
\end {cor}

\begin {rmk} After presenting these results at 
the Queen Dido conference on isoperimetry, we
learned that Laugesen and Morpurgo \cite{LM}
proved a series of very general results which 
includes the inequality of Theorem \ref{eigen-schwarz}. 
(See, in particular, Theorem 7, on page 80 of their
paper.) Our proof is quite different from that 
in \cite{LM} and may be of interest in its own right.
\end {rmk}

\begin {rmk}
One key point of \cite{BMMPR} is that their estimates 
involve the area 
(for instance) of the image of $f(r\D)$, rather than 
the area with multiplicity. 
We have left the corresponding question for the 
first Dirichlet eigenvalue 
of the Laplacian open. In this case, the first 
variation formula for the first 
eigenvalue of the image domain $f(r\D)$ is more 
complicated than the variation formula we 
have below, and when pulled back to $r\D$ will 
involve an integral over a proper subset of the 
boundary circle. 
\end {rmk}

A key step in the proof of this eigenvalue estimate 
is to rewrite a result of Payne and 
Rayner \cite{PR} as an isoperimetric-type inequality 
for the first eigenfunction. 
\begin{thmA} \label {eigen-isop}
Let $D$ be a bounded planar region with Lipschitz 
boundary $\partial D$, and let 
$\phi$ be the first eigenfunction of the Laplacian with Dirichlet boundary conditions. 
Then 
\beq\label{1.2}
\left ( \int_{\del D} |\nabla \phi| \right )^2 \geq 
4\pi \int_D |\nabla \phi|^2,
\eeq
with equality if and only if $D$ is a disk. 
\end{thmA}

The inequality \eqref{1.2} is in fact the 
isoperimetric inequality $L^2 \geq 4\pi A$ for the 
domain $D$, where one measures length 
and area with respect to the (singular) conformal 
metric $ds^2 = |\nabla \phi|^2 |dw|^2$. We discuss 
some properties of this metric below, in 
Section~\ref{eigen-sec}, and the equality in the 
case of the disk in Section~\ref{bessel-disk-sec}. 

The rest of this paper proceeds as follows. 
We prove Theorem~\ref{eigen-schwarz} in 
Section~\ref{eigen-sec} by writing out the first 
variation of the eigenvalue under a domain perturbation 
and reducing our problem to the isoperimetric 
inequality in Theorem~A. We examine the equality case 
of the isoperimetic inequality, that of a disk, 
in Section \ref{bessel-disk-sec}. 


\bigskip\noindent{\sc Acknowledgements:} We first learned about these variations 
on Schwarz's Lemma from Pietro Poggi-Corradini during a 
\textsl{Summer School in Conformal Geometry, Potential Theory, and Applications\/} 
at NUI Maynooth in June 2009. We would like to thank Poggi-Corradini 
for interesting discussions on the subject and the organizers of the conference for providing a stimulating venue for these discussions. We would also like to thank Michiel van den Berg 
for many enlightening conversations, and telling 
us of \cite{PR}. Finally, we would like to thank 
Rick Laugesen for pointing out \cite{LM}. 

\section{A Schwarz Lemma for the first eigenvalue}
\label{eigen-sec}

Let $f$ be an analytic function in the unit disk $\D$. 
Let $D_r = f(r\D)$ for $0<r<1$, and let $\lambda(r) 
= \lambda(D_r)$. 
The first Dirichlet eigenvalue for a disk is known to 
be  $\lambda(r\D) = j_0^2/r^2$ where $j_0$ is the 
first positive zero of the Bessel function $J_0$ of 
index zero. We then have 
\beq\label{2.1}
\Phi(r) = \frac{\lambda\big( f(r\D) \big)}
{\lambda(r\D)} = \frac{1}{j_0^2}\, r^2 
\lambda(r).
\eeq
Taking a derivative, we see that
\beq\label{2.2}
\frac{d\Phi}{dr} = \frac{1}{j_0^2}
\,\left [ 2r\lambda(r) + r^2 \frac{d\lambda}{dr}
\right ],
\eeq
so $\Phi$ is a decreasing function of $r$  precisely if 
\begin{equation} \label{decreasing-eigen-a}
\frac{2}{r}\,\lambda (r) \leq - \frac{d\lambda}{dr}. 
\end {equation} 

A classical theorem of Hadamard \cite{Had} computes 
the first variation of the 
eigenvalue as follows (see also \cite{Sch, E-SI, PS}). 

Let $\Omega_0$ be a domain. Let  $\zeta(t,x)$ be a 
flow on $\Omega_0$ associated 
with the variation field $\chi = \chi(t,x)$ in the 
time interval $(-t_0,t_0)$, in that,
\begin{align}
\frac{\partial \zeta}{\partial t} (t,x) & = 
\chi \big(\zeta(t,x)\big)\label{2.4}\\
\zeta(0,p) & = p, \quad p \in \Omega_0.
\label{2.5}
\end{align}
Let $\Omega_t$ be the domain $\zeta(t,\Omega_0)$, 
let $\lambda(t)$ be the first Dirichlet
eigenvalue for the Laplacian in $\Omega_t$, and let 
$\phi(t,x)$, $x \in \Omega_t$, 
be the associated eigenfunction normalised so 
that $\int_{\Omega_t} \phi^2  =  1$. 
Let $\eta$ denote the outward normal  and $d\sigma$ 
denote arc-length measure for $\del\Omega_t$.  
For the reader's convenience, we include a proof 
of the following formula for the time derivative 
of the eigenvalue, which draws heavily on the treatment 
in \cite{PS}.   
We take all boundaries and variation fields to be 
$C^\infty$, even though the variation formula holds 
with less regularity. In the calculation below 
we denote differentiation with respect to the parameter 
$t$ with a dot. 
\begin{lemma}\label{lemma3}
\begin {equation} \label {2.5a}
\dot \lambda(0)  = - \int_{\del \Omega_0} 
\langle \chi, \eta\rangle \left ( \frac
{\del \phi}{\del \eta} \right )^2 d\sigma . 
\end {equation} 
\end{lemma}

\begin {rmk} Because the first eigenvalue is 
simple, the function $\lambda(t)$ is differentiable. 
The higher eigenvalues $\lambda_k(t)$, for $k>1$, 
may not be differentiable functions of $t$, but 
both one-sided derivatives will exist. See the 
discussion in Sections 2 and 3 of \cite{E-SI} for 
more information. 
\end {rmk}

\begin {proof} First we compute the time derivative 
of the boundary terms 
of the normalized first eigenfunction $\phi$. Taking 
a derivative of the condition 
\[
\phi\big(t,\zeta(t,p)\big) = 0,\  p \in 
\del \Omega_ 0
\]
with respect to $t$ and using \eqref{2.4}, we 
obtain
\[
\dot\phi\big(t,\zeta(t,p)\big) + 
\big\langle \nabla\phi\big(t, \zeta(t,p)\big),
\chi(p) \big\rangle = 0.
\]
Here and later, the gradient refers only to the 
spatial derivative. Set $t=0$ and use the fact 
that $\phi$ is constant along $\del \Omega_t$ to 
obtain 
\begin {equation} \label {first-var-a} 
\dot \phi(0,p)  
= - \big\langle \nabla \phi (0,p), \chi(p) 
\big\rangle = -\Big\langle \left.\frac{\del \phi} 
{\del \eta} \right|_{(0,p)} \eta(p), \chi(p)
\Big\rangle, \quad p \in \del \Omega_0.
\end {equation} 

Next we take the derivative of the eigenfunction 
equation 
\begin {equation} \label {first-var-b} 
\Delta \phi\big(t,\zeta(t,p)\big) + \lambda(t) \phi\big(t,\zeta(t,p)\big) = 0 
\end {equation} 
with respect to $t$ and evaluate at $t=0$. This leads to
\begin {eqnarray*}
0 & = & \Delta \left [ \dot \phi + \langle \nabla 
\phi, \chi \rangle \right] + \lambda(t) 
\left [ \dot \phi + \langle \nabla \phi, 
\chi\rangle \right ] + \dot \lambda(t) \phi\\
& = & \Delta \dot\phi + \langle \nabla \Delta \phi, 
\chi \rangle  +  \lambda(t) 
\dot\phi + \lambda(t) \langle \nabla \phi, 
\chi\rangle + \dot \lambda(t) \phi \\
& = & \Delta \dot\phi + \lambda(t) \phi_t + 
\dot\lambda(t) \phi .
\end {eqnarray*}
Setting $t=0$ and rearranging yields 
\begin {equation} \label{first-var-c} 
\Delta \left. \dot\phi \right|_{t=0} + \lambda(0) 
\left. \dot\phi \right|_{t=0} 
= -\dot\lambda(0)  \phi\big\vert_{t=0} \quad 
\text{ in } \Omega_0.
\end {equation} 

We multiply  (\ref{first-var-b}), with $t=0$,  
by $\dot\phi\big\vert_{t=0}$ and 
multiply (\ref{first-var-c}) by $\phi$, subtract and 
obtain 
\begin {equation} \label {first-var-d}
\dot\lambda(0) \phi^2(0,p) = \dot\phi(0,p) 
\Delta \phi(0,p) - \phi(0,p) 
\Delta \dot\phi(0,p), \quad p \in \Omega_0.
\end {equation} 

Integrate \eqref{first-var-d} over $\Omega_0$ and 
use the fact that $\int_{\Omega_t} \phi^2  = 1$ 
to obtain 
\begin {eqnarray*} 
\dot\lambda(0) & = & \int_{\Omega_0} \dot\phi \Delta \phi 
- \phi \Delta \dot\phi\\
& = & \int_{\del \Omega_0} \dot\phi 
\frac{\del \phi}{\del \eta} - 
\int_{\Omega_0} \big\langle \nabla \phi, 
\nabla \dot\phi \big \rangle + 
\int_{\Omega_0} \big\langle \nabla \phi, 
\nabla \dot\phi \big\rangle - \int_{\del \Omega_0} 
\phi \frac{\del \dot\phi}{\del \eta} \\
& = &\int_{\del \Omega_0} \dot\phi 
\frac{\del \phi}{\del \eta}\\
& = & - \int_{\del \Omega_0} \frac{\del \phi}
{\del \eta} \langle \nabla \phi, \chi \rangle \\
& = & - \int_{\del \Omega_0} \langle \chi, 
\eta\rangle \left ( \frac{\del \phi}{\del \eta} 
\right )^2, 
\end {eqnarray*} 
which is equation (\ref{eigen-deriv1}) as claimed. In the second equality above we 
integrated by parts, in the next to last we used \eqref{first-var-a}, and at the last step 
we used the fact that $\phi$ is constant on $\del \Omega_0$ (and hence $\nabla \phi = 
\frac{\del \phi}{\del \eta} \eta$ there). 
\end {proof} 

We adapt this formula to our particular case. 
\begin{lemma}Let $f$ be a conformal mapping of  
the unit disk $\D$ with $f(0)=0$. 
Let $\lambda(r)$ be the eigenvalue of the domain 
$D_r = f(r\D)$ with 
eigenfunction $\phi_r$ in $L^2(D_r)$. Let 
\begin{equation}\label{2.6}
\psi_r(z) = \phi_r\big( f(z)\big), \quad z \in r\D.
\end{equation}
Then
\begin {equation} \label {eigen-deriv1}
\frac{d\lambda}{dr} = - r\int_0^{2\pi} 
\vert(\nabla \psi_r)(re^{i\theta})\vert^2\, d\theta.
\end {equation} 
\end{lemma}

\begin {rmk} 
The function $\psi$ satisfies the equation 
$$\Delta \psi + \lambda |f'|^2 \psi = 0,$$
and so $\psi$ is the first Dirichlet eigenfunction of the 
Laplacian on the conformal disk $(\D, |f'|^2 |dz|^2)$. 
\end {rmk}

\begin{proof}
For a fixed $r$ in $(0,1)$, we set 
\begin{equation}\label{2.7}
\zeta(t,p) = f\big( (1+t/r)f^{-1}(p)\big), \quad p 
\in D_r.
\end{equation}
Then, $\zeta(0,p) = p$, $\zeta(0,D_r) = \Omega_0 = D_r$, $\zeta(t,D_r) = \Omega_t = 
D_{r+t}$, and 
\[
\frac{\del \zeta}{\del t}(t,p) = f'\big( (1+t/r)f^{-1}(p)
\big) \frac{1}{r}f^{-1}(p). 
\]
It follows that \eqref{2.4} holds with
\begin{equation}\label{2.8}
\chi (\zeta) = \frac{1}{r+t} f^{-1}(\zeta) 
f'\big(f^{-1}(\zeta)\big) .
\end{equation}
The unit normal vector to the boundary of $D_r$ at 
$\zeta$ is  
\begin{equation}\label{2.9}
\eta(\zeta) = \frac{f^{-1}(\zeta)}{r}\, 
\frac{f'\big(f^{-1}(\zeta)\big)}{
\big\vert  f'\big(f^{-1}(\zeta)\big) \big\vert}. 
\end{equation}
Thus, \eqref{2.8} with $t=0$ and \eqref{2.9} show that 
$\chi(\zeta) = \big\vert  f'\big(f^{-1}(\zeta)\big) 
\big\vert \eta(\zeta)$, $\zeta \in \del D_r$, 
so that 
\[
\langle \chi, \eta\rangle = \big\vert  
f'\big(f^{-1}(\zeta)\big) \big\vert, \quad
\zeta \in \del D_r.
\]
The gradient of the function $\psi_r$ given by 
\eqref{2.6} is $ \vert \nabla \psi_r(z)\vert = 
\vert \nabla \phi_r\big( f(z) \big) \vert\, \vert 
f'(z)\vert$. This, and Lemma~\ref{lemma3}, lead to 
\begin{align*}
\dot\lambda(r) &  = - \int_{\del D_r}  \big\vert  
f'\big(f^{-1}(\zeta)\big) \big\vert\,
\vert \nabla \phi_r(\zeta) \vert^2\, \vert 
d\zeta\vert\\
& = -\int_{C(0,r)} \vert f'(z)\vert \left( 
\frac{\vert \nabla \psi_r(z)\vert}{\vert f'(z)
\vert}\right)^2 \, \vert f'(z)\vert \, \vert dz\vert\\
& =  -\int_{C(0,r)} \vert \nabla \psi_r(z)\vert^2 
\, \vert dz\vert,
\end{align*} 
where $C(0,r)$ denotes the circle centre 0 and 
radius $r$, which is \eqref{eigen-deriv1}. 
\end{proof}

Next we verify Theorem~A. 
Payne and Rayner write the inequality in the form
\begin{equation}\label{1.4}
\left( \int_D \phi \right)^2 \geq 
\frac{4\pi}{\lambda} \int_D \phi^2 
\end{equation}
with equality if and only if $D$ is a disk. 
If we denote by $\eta$ the unit outward normal to 
the boundary of $D$,  
\[
\int_{\del D} \vert \nabla\phi\vert = \int_{\del D}
\left( - \frac{\del\phi}{\del\eta}\right)
 = \int_{D}\left( - \Delta\phi\right)  =  \lambda 
\int_D\phi,
\]
where the first equality comes from the fact that 
$\phi$ is constant on the boundary
of $D$, the second from Green's theorem, and the 
third from the eigenfunction 
equation $\Delta \phi + \lambda \phi = 0$. Since 
$\phi$ minimises the Rayleigh quotient,
\[
\int_D \vert \nabla \phi \vert^2 = \lambda \int_D \phi^2.
\]
Hence,  \eqref{1.2} and \eqref{1.4} are equivalent.

\begin{proof}[Proof of Theorem~\ref{eigen-schwarz}]
By Theorem~A,
\[
\left( \int_{C(0,r)} \big\vert \nabla \psi_r(z)\big
\vert\,\vert dz\vert \right)^2
 = \left( \int_{\del D_r)} \vert \nabla \phi_r\vert 
\right)^2 \geq 4\pi \int_{D_r} \vert \nabla 
\phi_r\vert^2 =  4\pi \int_{r\D} \vert \nabla 
\psi_r\vert^2.
\]
Hence, since $\Vert \phi_r \Vert_{L^2(D_r)}=1$, 
\begin{align*}
\frac{2}{r}\,\lambda (r)  = \frac{2}{r} \int_{D_r} 
\vert \nabla \phi_r\vert^2 
& = \frac{2}{r} \int_{r\D} \vert \nabla 
\psi_r\vert^2\\
& \leq \frac{r}{2\pi} \left( \int_0^{2\pi} 
\vert \nabla \psi_r(re^{i\theta})\vert
\,d\theta\right)^2\\
& \leq r \int_0^{2\pi} \vert \nabla \psi_r
(re^{i\theta})\vert^2\,d\theta\\& = -\lambda'(r),
\end{align*}
by \eqref{eigen-deriv1}. This proves \eqref{decreasing-eigen-a} and hence 
Theorem~\ref{eigen-schwarz}. 
\end{proof}

\begin {rmk} 
The metric $|\nabla \phi|^2 |dw|^2$ on $D_r = 
f(r\D)$ is a singular metric, with singularities 
at the critical points of the eigenfunction $\phi$. 
However, $\phi$ solves a second order, linear, 
elliptic equation on a bounded domain with $C^\infty$
boundary, so it only has finitely many critical 
points, none of which are degenerate. We 
conclude that the metric $|\nabla \phi|^2 
|dw|^2$ in $D_r$, or, equivalently, $|\nabla \psi|^2
|dz|^2$ in $r\D$, has only finitely many singular 
points, where the metric vanishes only to zeroth 
order. 
\end {rmk}

\begin {rmk}
The isoperimetric inequality $L^2 \geq 4\pi A$, 
with $L = \int_{\del D} |\nabla \phi| |dw|$ and 
$A = \int_D |\nabla \phi|^2 |dw|^2$, seems quite 
general. Using the Riemann mapping theorem, 
it holds for any simply connected domain $D 
\subset \C$. We also find it remarkable that 
the constant in this isoperimetric inequality 
is the same one as in the classical isoperimetric 
inequality. 
\end {rmk} 

\begin{proof}[Proof of Corollary~\ref{eigen-schwarz2}]
The function $f: r\D \rightarrow \Sigma_r$ is 
a conformal map away from its critical points, and 
so the first variation formula 
\eqref{eigen-deriv1} holds so long as $f$ does not 
have a critical point of length $r$. For any $r_0 
\in (0,1)$ there are only finitely many values 
$\hat r < r_0$ such that $f$ has a critical point of 
length equal to $\hat r$, and the variation 
formula \eqref{eigen-deriv1} is valid away from 
these values $\hat r$. Thus, we can integrate 
the inequality \eqref{decreasing-eigen-a} to see 
that $\Phi(r) = (r^2/j_0^2) \lambda(f(r\D))$ is 
decreasing for $0<r<r_0$. Moreover, if there 
exist $r_1 < r_2$ such that $\Phi(r_1) = \Phi(r_2)$
then $f$ is linear on the annulus $r_1 < |z| < r_2$; 
this combined with the fact that $f$ is analytic 
on the disk $\D$ implies $f$ is linear on the 
whole disk. 
\end {proof} 

\section {The Bessel disk: the equality 
case of the isoperimetric inequality}
\label {bessel-disk-sec}

In any given inequality, the case of 
equality is always important, and often 
sheds light on other problems. The 
equality case of \eqref{eigen-isop} 
occurs when $f$ is linear, in which case 
the image domains $f(r\D)$ are disks for all 
$r \in (0,1)$. In this case, the eigenfunctions 
$\phi$ and $\psi$ agree up to scaling factors, 
and we write $\phi(z) = J_0 (j_0 |z|)$, 
where $J_0$ is the Bessel function with 
index zero and $j_0$ is its first 
positive root. 

\begin {defn} We call the unit disc $\D$ 
equipped with the conformal metric 
$$ds = J_1(j_0|z|) |dz|$$ 
the Bessel disc.
\end {defn} 
The following lemma is an immediate 
consequence of \eqref{eigen-isop}.
\begin {lemma} In the Bessel disc, $L^2 
= 4\pi A$. 
\end {lemma} 

We conclude this section by exploring 
the geometry of the Bessel disk. 
It is convenient to recall the formula 
for the Gauss curvature of a conformal 
metric. In general the metric $ds = \rho 
|dz|$ has Gauss curvature 
$$K  = -\frac{1}{\rho^2} \Delta 
(\log \rho) .$$
In particular, negative curvature is 
equivalent to $\log \rho$ being subharmonic. 

Observe that the Bessel disk has positive 
curvature, which blows up logarithmically 
at the origin. We include plots 
of the curvature and the Gauss-Bonnet 
integrand for the reader's enlightenment.

\begin {center} 
\begin {figure}
\includegraphics [width = 4in]{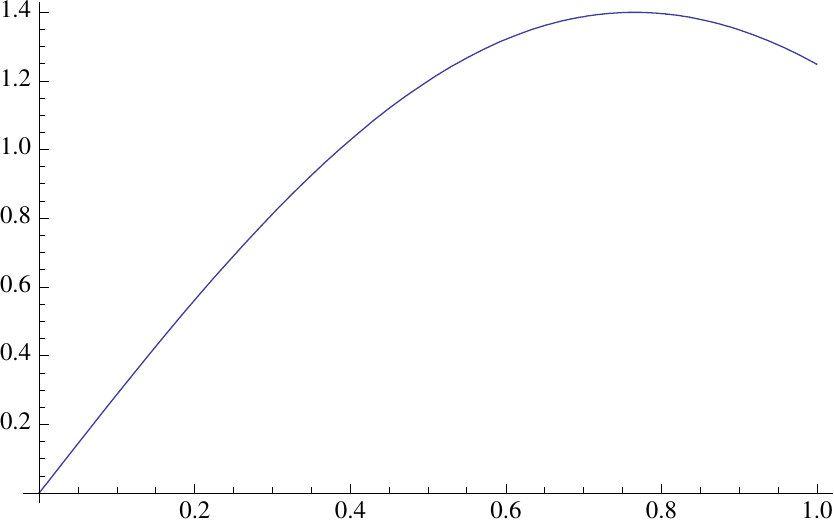}
\caption  {This is a plot of the 
conformal factor $\rho(r) = 
j_0 J_1(j_0 r)$ for the Bessel disk.}
\end {figure} 
\end {center} 

\begin {center} 
\begin {figure}
\includegraphics [width = 4in]{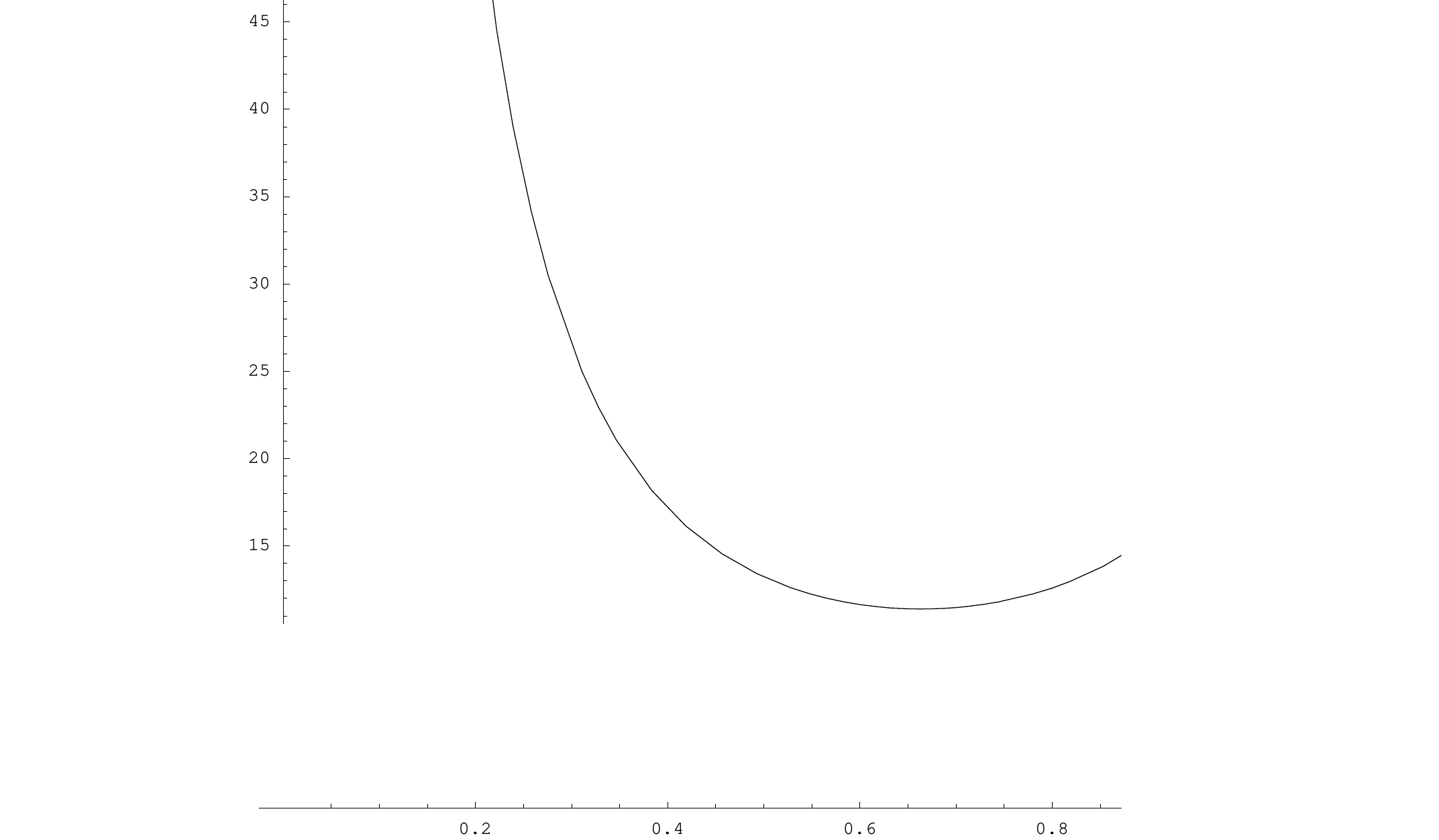}
\caption  {This is a plot of the 
curvature $K = 
-\rho^{-2} \Delta \log(\rho)$ for the Bessel disk.}
\end {figure} 
\end {center} 

\begin {center} 
\begin {figure}
\includegraphics [width = 4in]{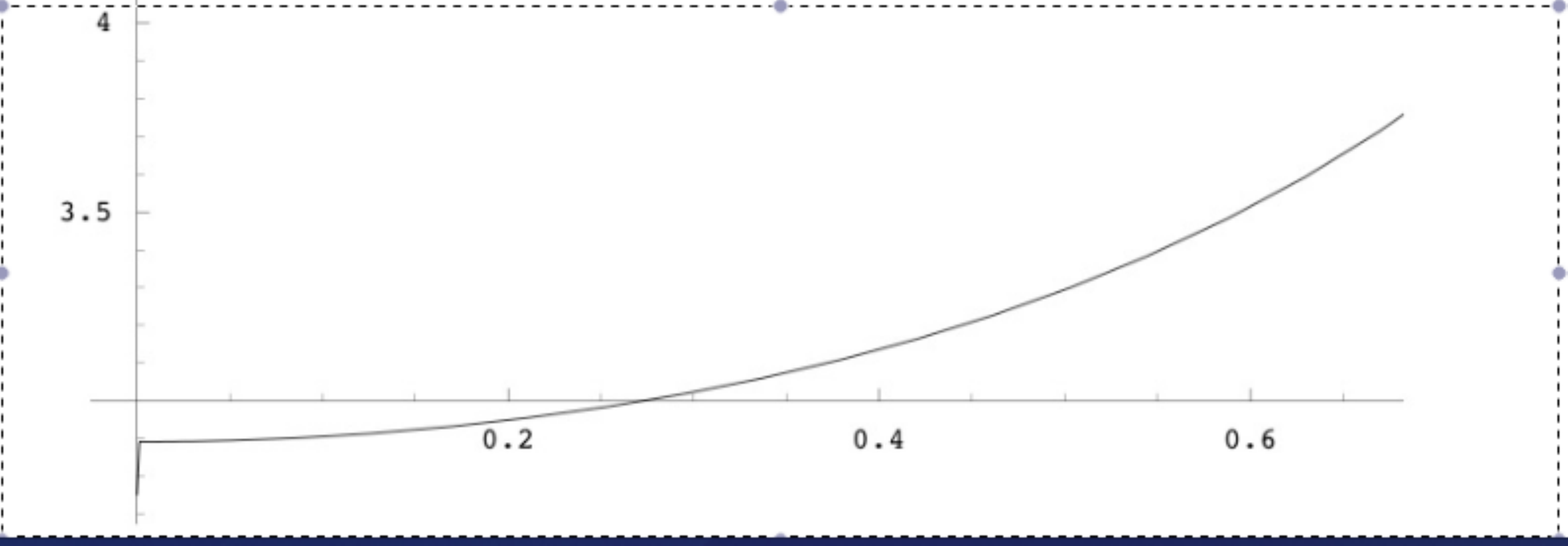}
\caption  {This is a plot of the 
Gauss-Bonnet integrand $K dA = -\Delta \log (\rho)|dz|^2 $ for the Bessel disk.}
\end {figure} 
\end {center} 

\begin {lemma} The total curvature 
of the Bessel disc is $4\pi$. 
\end {lemma} 

\begin {proof} Let $\rho (z) = 
J_1(j_0 |z|)$, so that 
\begin {eqnarray*} 
\int_{\D} K dA & = & -\int_{\D} \frac
{\Delta \log \rho}{\rho^2} \rho^2 |dz|^2 
= -2\pi \int_0^1 [(\log \rho)''(r) 
+ r^{-1} (\log \rho)'(r) ] rdr \\
& = & -2\pi \int_0^1 r(\log \rho)''(r) 
+ (\log \rho)'(r) dr = -2\pi \int_0^1 
\frac{d}{dr} (r (\log \rho)') dr \\ 
& = & 2\pi \left. \left ( \frac{r
\rho'(r)}{\rho(r)} \right ) \right|_0^1
= 2\pi \left ( \lim_{r\rightarrow 0} 
\frac{r j_0 J_1'(j_0 r)}{J_1(j_0 r)} 
- \frac{j_0 J_1'(j_0)}{J_1(j_0)} \right ) \\ 
& = & 4\pi.
\end {eqnarray*} 
Here we have used Bessel identities to  
show $j_0 J_1'(j_0) = - J_1(j_0).$ 
\end {proof}

Looking closely at this computation, we
see that the Gauss-Bonnet integrand is 
an exact derivative, and so there are 
two terms which contribute to the total 
curvature: a boundary term and an interior 
term at the critical point of the first 
eigenfunction. For any bounded domain 
$D$ with Lipschitz boundary, the local 
behavior of its first eigenfunction near a 
critical point will be that of the Bessel 
function at the origin of the disk, at least 
to first order. Thus, the computation above shows 
that any critical point of the first eigenfunction 
will contribute $2\pi$ to the total curvature of 
$(D, |\nabla \phi|^2 |dz|^2)$, where $\phi$ is the first 
eigenfunction of $D$. It therefore seems natural 
to conjecture that, for instance, 
the total curvature of $(D, |\nabla \phi|^2 |dz|^2)$ is 
exactly $4\pi$ for any convex domain $D$. 

We contrast the Bessel disk with the 
isoperimetric inequalities of Topping 
\cite{Top} and Morgan-Hutchings-Howard 
\cite{MHH}. They prove that a rotationally 
symmetric metric with a monotone curvature 
function will achieve equality in each of 
their inequalities. On the other hand, one can 
verify the following properties of the 
Bessel disk by explicit computation. First, 
it is a rotationally symmetric metric, 
which realizes equality in both the 
isoperimetric inequalities of \cite{Top} and 
\cite{MHH}. Second, the curvature is {\emph not} 
monotone. It remains an interesting open question 
to characterize which metrics achieve equality 
in the isoperimetric  inequalities of \cite{Top} 
and \cite {MHH}.

\end{document}